\documentstyle[11pt,righttag,intlim]{amsart}
\title{A commutator method for computation of heat invariants}
\author{Iosif Polterovich}
\address{Department of Theoretical Mathematics, The Weizmann Institute
of Science, Rehovot 76100, Israel}
\email{iossif@@wisdom.weizmann.ac.il}
\date{}
\def \phi{\varphi}
\def \epsilon{\varepsilon}
\numberwithin{equation}{section}
\theoremstyle{definition}
\newtheorem{definition}[equation]{Definition}
\theoremstyle{plain}
\newtheorem{lemma}[equation]{Lemma}
\newtheorem{theorem}[equation]{Theorem}

\begin{document}
\maketitle
\begin{abstract}
We introduce a new method  for computing heat invariants $a_n(x)$ of a 
$2$-dimensional Riemannian manifold  based on a commutator formula derived
by S.Agmon and  Y.Kannai. Two explicit expressions for $a_n(x)$ 
are presented. The first
one depends on the choice of a certain coordinate system;
the second involves only invariant terms but has some restrictions
on its validity, though in  a ``generic'' case it is well--defined. 
\end{abstract}

\section{Introduction and main results}
Let $M$ be a $d$-dimensional compact Riemannian manifold 
without boundary with a  metric $\{g_{ij}\}$. 
Let $\Delta$ be the Laplace--Beltrami operator (or simply the Laplacian)
on $M$.
In the coordinate chart  $\{u_i\}$ it is written in the form
\begin{equation}
\Delta f = - \frac{1}{\sqrt{g}}\sum_{i,j=1}^d 
\frac{\partial(\sqrt{g} g^{ij}(\partial f/\partial u_i))}{\partial u^j},
\end{equation}
where $g=$det$\{g_{ij}\}$, and $\{g^{ij}\}$ denotes the inverse of the 
matrix $\{g_{ij}\}$.

The Laplacian is a self-adjoint elliptic operator acting on
smooth functions $f: M \to {\Bbb R}$. 
Let 
$\{ \lambda_i, \phi_i \}_{i=1}^\infty$ be the spectral decomposition 
of $\Delta$
(considered as an operator in $L^2(M)$) 
into a complete orthonormal basis of
eigenfunctions $\phi_i$ and eigenvalues $\lambda_i$, $0 \le \lambda_1 \le
\lambda_2 ...$.
Consider the heat operator $e^{-t\Delta}$ which is the solution of the 
heat equation $(\partial /\partial t+\Delta)f=0$. For $t$ positive,
it is an infinitely
smoothing operator from $L^2(M)$ into itself and has a  
kernel function 
$${\mbox K}(t,x_1,x_2)=\sum e^{-\lambda_i t} \phi_i(x_1)\otimes \phi_i(x_2),$$
analytic in $t$
and $C^{\infty}$ in 
$x_1$ and $x_2$, such that $$(e^{-t\Delta}f)(x_1)=\int\limits_M 
{\mbox K}(t,x_1,x_2)f(x_2)dx_2.$$ 
The following asymptotic expansion holds for the trace of the
heat kernel:
$${\mbox K}(t,x,x) \sim \sum_{n=0}^{\infty} a_n(x)t^{n-\frac{d}{2}}
$$
as $t \to 0^+$ ([G1]).
The coefficients $a_n(x)$ are called {\it heat invariants} of the manifold $M$;
they depend only on the germs of the metric. 
Moreover, it is known that the heat invariants are 
homogeneous polynomials of degree $2n$ in the derivatives of 
$\{g^{ij}\}$ at the point $x$ ([G2]).

Integrating $a_n(x)$ over
the manifold one gets the coefficients of the expansion of the trace 
of the heat operator:  
$$\sum_i e^{-t\lambda_i} \sim  \sum_{n=0}^{\infty}\left(\int\limits_M 
a_n(x,\Delta)\rho(x)dx \right) t^{n-\frac{d}{2}} \sim  \sum_{n=0}^{\infty}a_n
t^{n-\frac{d}{2}}, $$
where $\rho(x)dx$ is the volume form on the manifold $M$.

Heat invariants are metric invariants of the given manifold $M$ and
contain a lot of information about its geometry and topology 
(see [Be], [G2], [G3], etc.).
For example, $a_0$ determines the volume of $M$, $a_1$ its scalar 
curvature (and hence the Euler characteristic in dimension 2). 
Higher coefficients are of a considerable interest to
physicists since they are connected with many notions 
of quantum gravity (see [F]).

The structure of the heat invariants becomes more and more complicated as 
$n$ grows; compare the formula for $a_3$ 
([Sa], 1971) with the formulas for $a_4$ 
([ABC], 1989; [Av], 1990)  and $a_5$ ([vdV], 1997).
For $n>5$ no explicit expression for the heat invariants is known, even though
there are many algorithms to derive them --- the computational complexity
is so high that even with the help of modern computers it is impossible
to get the result.

In the present paper we introduce a new method for computing heat
invariants. We call it the commutator  method since it is 
substantially based on a commutator formula by S. Agmon and Y. Kannai ([AK])
which deals with asymptotic expansions of
the resolvent kernels of elliptic operators (let us note that the 
asymptotics of the relevant kernels associated with a 
self--adjoint realization of the elliptic operator depend only on the 
local behavior of its coefficients).

In section 2 we obtain the following concise reformulation of 
the Agmon--Kannai formula:  
\begin{theorem}
Let $H$ be a a self--adjoint elliptic differential 
operator of order $s$ on a manifold $M$ of dimension $d<s$, 
and $H_0$ be the operator obtained by
freezing the coefficients of the principal part $H'$ of the operator
$H$ at some point $x\in M$: $H_0=H_x'$.  
Denote by $R_\lambda(x_1,x_2)$ the kernel of the resolvent
$R_\lambda=(H-\lambda)^{-1}$, and by $F_{\lambda}(x_1,x_2)$ --- the 
kernel of $F_{\lambda}=(H_0-\lambda)^{-1}$.
Then for $\lambda \to \infty$ 
there is a following asymptotic representation on the diagonal:
\begin{equation}
R_\lambda(x,x) \sim \rho(x)^{-1} \sum_{m=0}^\infty X_m F_\lambda^{m+1}(x,x),
\end{equation}
and the operators $X_m$ are defined inductively by:
$$
X_0=I; \,\,\, X_m=X_{m-1}H_0-HX_{m-1}.
$$
\end{theorem}

The commutator  method may be applied to various problems 
of ``heat kernel''--type 
(one may consider other operators instead of Laplacian,
Laplacian on forms, etc. --- see also [Kan] for the Agmon--Kannai 
formula in case of matrix operators).
In this  paper we are focusing on the simplest non--trivial question ---
what are the heat invariants of a 2-dimensional Riemannian manifold?

An advantage of the commutator  method is that it 
yields very explicit formulas for {\it all} 
2-dimensional heat invariants (cf. [Xu] where heat invariants 
of manifolds of arbitrary dimension are considered).
They are given by the following 
\begin{theorem}
Let $(u,v)$ be local coordinates in a neighborhood of 
the point $x=(0,0)$ on the Riemannian $2$-manifold $M$,
in which the metric is conformal: $$ds^2=\rho(u,v)(du^2+dv^2).$$
Then the heat invariants $a_n(x)$ are given by 

\begin{equation}
a_n(x)=\sum_{m=n+1}^{4n}\sum_{k=n+1}^m \sum_{s=0}^{k-n}C_{nksm}
\rho(u,v)^{k-n}
\Delta^k\left(u^{2k-2n-2s}v^{2s}\right)|_{u=0,v=0},
\end{equation}
where $C_{nksm}$ are constants given by {\em (3.12)}.
\end{theorem}
The proof of this theorem is presented in section 3.
We also perform a test of the formula (1.5) which is described in the
Appendix.

Let us note that local coordinates in which the metric is
conformal always exist (see [DNF]). 
We also observe that the derivatives of $\rho(u,v)$ appear in (1.5) 
through the formula (1.1) for the coefficients of the Laplacian.

The expression (1.5) depends on a coordinate representation but
it is possible to make it invariant rewriting it in special
curvature coordinates which are introduced in section 4.
The resulting formula is no longer valid in general but
is true ``generically''.

\medskip

\noindent {\it Acknowledgments.} This paper is a part of my
Ph.D. research conducted under the supervision of Professor Y. Kannai.
I am very  grateful to Professor Kannai for his constant help  and 
encouragement.  I would also like to thank Professor L. Polterovich and 
Professor M. Solomyak for many fruitful discussions, and 
my colleagues J. Greenstein and  A. Grigoriev for helpful advice.
\medskip

\section{Multiple commutators}
We recall the notion of a multiple commutator of two operators, introduced
in ([AK]), which plays a fundamental role in the sequel.
\begin{definition}
Let $J$ be a  finite vector with nonnegative integer components
and let $A$, $B$ be linear operators on some linear space $M$.
The {\em multiple commutator} $[B,A;J]$  is inductively defined for all
such vectors $J$ in the following way:
$$1.\quad [B,A;0]=B$$
$$2.\quad [B,A;j]=[[B,A;j-1],A]$$
$$3.\quad [B,A;J\cup j]=[B[B,A,J],A;j],$$
where $J\cup j$ denotes the vector obtained by adding the component $j$ 
to the vector $J$ to the right.
\end{definition} 
We introduce a filtration on the space of all vectors  
with non-negative integer entries.
Denote by $V_m=\{J=(j_1,..,j_r):|J|+r=m\}$, where $|J|=j_1+..+j_r$.
Let
\begin{equation}
X_m= \sum_{J \in V_m}[B,A;J], \quad m \ge 1.
\end{equation}

\begin{theorem}
The operators $X_m$ satisfy the following recurrent relation,
\begin{equation}
X_1=B; \quad X_m=BX_{m-1}+[X_{m-1},A],
\end{equation}
and are given explicitly by 
\begin{equation}
X_m=\sum_{k=0}^m (-1)^k \binom{m}{k} (A-B)^k A^{m-k}.
\end{equation}
\end{theorem}
\noindent {\bf Proof.}
Let us  proceed by induction.
Indeed, $X_1=B$ since $V_1$ consists of a single vector $(0)$.
In order to pass from $X_{m-1}$ to $X_m$ we do the following.
Let $V_m^0$ denote the set of all vectors in $V_m$ whose last entry
is $0$. Consider a mapping $p:V_m \to V_{m-1}$ defined as follows:
for $J=(j_1,..,j_{r-1},0) \in V_m^0$ we have $p(J)= (j_1,..,j_{r-1})$
and for $J=(j_1,..,j_r) \in V_m\setminus V_m^0$ we have $p(J)=(j_1,..,j_r-1)$.
It is clear that the mapping $p$  is well--defined. Moreover, it is 
a surjection and every element of $V_{m-1}$ has exactly two preimages ---
one in $V_m^0$ and one in its complement.
Applying the definition of a multiple commutator we prove 
the inductive step --- $BX_{m-1}$ is the contribution of vectors 
from $V_m^0$ and $[X_{m-1},A]$ is the contribution of vectors belonging 
to $ V_m\setminus V_m^0$. This completes the proof of the first part
of the theorem. 
For obtaining the closed formula, let us rewrite the recurrent relations
in the following way:
$X_1=A-(A-B)=B$; $X_m=X_{m-1}A-(A-B)X_{m-1}$.
Proceeding by induction and recalling that 
$\binom{m}{k}+\binom{m}{k-1}=\binom{m+1}{k}$ we get (2.5) \qed

\bigskip

Let us present the original Agmon--Kannai formula (see [AK]):

\begin{theorem}
In the notations and conditions of the Theorem {\em 1.2} the following
asymptotic representation on the diagonal holds for the kernel
of the resolvent $R_{\lambda}$:
\begin{equation}
R_\lambda(x,x) \sim \rho(x)^{-1}(F_\lambda(x,x)+ 
\sum_{J}^\infty ([H_0-H,H_0;J]F^{|J|+r+1}(x,x))), 
\end{equation}
where the sum is taken over all vectors $J$ of length $\ge 1$ with 
nonnegative integer entries.
\end{theorem}
Now we may prove the Theorem 1.2.

\medskip

\noindent {\bf Proof of  the Theorem 1.2.}
Let us divide the above sum into sums over 
vectors $J \in V_m$, $m \ge 1$.
Then, setting $A=H_0$, $B=H_0-H$ and applying Theorem 2.3.
we get formula (1.3), which completes the proof. \qed

\section{The commutator method in dimension 2}
Theorem 1.2 is not valid for the Laplacian if $d \ge 2$.
A way to overcome this difficulty for $d=2$ is to consider
differences of resolvents (this was 
pointed out to me by  Y.Kannai). 
  
\begin{lemma}
Let $R_{\lambda}=(\Delta -\lambda)^{-1}$ be the resolvent of the Laplacian
on a $2$-manifold $M$.
Consider the difference of the resolvents
$R_{\lambda,2\lambda}=R_\lambda-R_{2\lambda}$.
This operator has a continuous kernel $R_{\lambda,2\lambda}(x_1,x_2)$
which has the following 
asymptotic representation on the diagonal as $\lambda \to \infty$
\begin{equation}
R_{\lambda,2\lambda}(x,x) \sim \sum_{n=0}^{\infty}b_n(x)
(-\lambda)^{-n},
\end{equation}
where  
\begin{equation}
\label{bn}
b_n(x)=(n-1)! \cdot a_n(x), \quad n>1.
\end{equation}
\end{lemma}

\noindent {\bf Proof.}
The difference of the resolvents $R_{\lambda,2\lambda}$ 
is a self--adjoint smoothing operator from 
$L^2(M)$ into the Sobolev space ${\mbox H}^4(M)$, 
and $\operatorname{dim} M < 4$, therefore
it has a continuous kernel (see [AK]).

Let us now prove the expansion (3.2). 

Set $z=-\lambda$, and let $Re z>0$.
We have (formally):
$$\int_0^\infty e^{-t(\Delta+z)}dt=\frac{1}{z+\Delta}$$
Therefore
$$
\frac{d}{dz}(\frac{1}{z+\Delta})=-\int_0^\infty te^{-t(\Delta+z)}dt
\sim -\int_0^\infty t\sum_{n=0}^\infty t^{-1+n}a_ne^{-zt}dt.
$$
Let us integrate both parts from $z$ to $2z$:
$$
\frac{1}{z+\Delta}-\frac{1}{2z+\Delta} \sim
\int_z^{2z}\int_0^\infty\sum_{n=0}^\infty a_n t^n e^{-zt}dt.
$$
For $n=0$ we get:
$$\int_z^{2z}a_0\int_0^\infty e^{-zt}dt dz=\int_z^{2z} \frac{a_0}{z} dz=
a_0 \log 2 = a_0'$$
For $n\ge 1$ we get:
$$
\int_z^{2z}\int_0^\infty t^n e^{-zt} dt dz= \Gamma(n+1)a_n\int_z^{2z}
\frac{dz}{z^{n+1}}=(n-1)! a_n (\frac{1}{z^n}-\frac{1}{(2z)^n}).
$$
This completes the proof of the lemma. \qed

\begin{lemma}
Let $R_\lambda(x_1,x_2)$ be the kernel of the resolvent
$R_\lambda=(\Delta-\lambda)^{-1}$
Then for $\lambda \to \infty$ 
the difference of the resolvents has the following asymptotic representation 
on the diagonal:
\begin{equation}
\label{sum}
R_{\lambda,2\lambda}(x,x) \sim \rho(x)^{-1} \sum_{m=0}^\infty X_m 
(F_\lambda^{m+1}(x,x)-F_{2\lambda}^{m+1}(x,x)),
\end{equation}
where 
\begin{equation}
X_m=\sum_{k=0}^m (-1)^k \binom mk \Delta^k\Delta_0^{m-k}. 
\end{equation}
\end{lemma}

\noindent {\bf Proof.}
Setting  $\Delta=H$ and applying  Theorem 1.2 and Lemma 3.1 
completes the proof.
Note that we can define $X_0=I$ without loss of consistence of
the recurrent formula (2.4.). \qed

Now we are ready to prove an explicit formula for the heat invariants
in two dimensions.

\medskip

\noindent {\bf Proof of the Theorem 1.4.}
Let  $(u,v)$ be conformal local coordinates in the neighborhood
of the point $x=(0,0)$ (we can assume that $x$ is the 
origin without loss of generality) on the $2$-manifold $M$.
As it was mentioned in the introduction this is always possible.
The metric in such coordinates has the form
\begin{equation}
ds^2=\rho(u,v)(du^2+dv^2),
\end{equation}
and hence the Laplacian is given by
\begin{equation}
\Delta= -\rho(u,v)^{-1}(\partial^2/\partial u^2+\partial^2/\partial v^2).
\end{equation}
Our aim is to find the coefficient $a_n(x)$. 
In order to do this we have to collect all
terms in the sum (\ref{sum})) containing $(-\lambda)^{-n}$. 
We use the following well--known formula (see [AK]):
\begin{equation}
\label{int}
\frac{\partial^{\alpha_1+\alpha_2}}{\partial u^{\alpha_1} 
\partial v^{\alpha_2}}
F_\lambda^{m+1}(x,x)=
(-\lambda)^{\frac{|\alpha|}{2}-m}(-1)^{\frac{|\alpha|}{2}}
\frac{1}{4\pi^2}\int\limits_{\Bbb R^2}
\frac{\xi_1^{\alpha_1}\xi_2^{\alpha_2}\,d\xi_1d\xi_2}
{(\Delta_x'(\xi_1,\xi_2)+1)^{m+1}}, 
\end{equation}
where $\Delta_x'(\xi_1,\xi_2)$ denotes the symbol of the operator
$\Delta_x'$,  and
$(-\lambda)^{\frac{|\alpha|}{2}-m}$ is the analytic branch 
of the power which is positive on the negative axis.

By formula (\ref{int}) we get $|\alpha|=2(|J|+r-n)=2(m-n)$.
Since the integral (\ref{int}) vanishes unless both $\alpha_1$
and $\alpha_2$ are even we may put $\alpha_1=2m-2n-2p$, $\alpha_2=2p$,
where $0\le p\le 2(m-n)$.
In polar coordinates this integral
transforms into the product of two integrals:
\begin{multline*}
\int\limits_{\Bbb R^2}
\frac{\xi_1^{2m-2n-2p}\xi_2^{2p}}
{(\Delta_x'(\xi_1,\xi_2)+1)^{m+1}} d\xi_1d\xi_2=\\
=\frac{1}{2}\int\limits_0^{2\pi}\cos^{2m-2n-2p}\phi \sin^{2p}\phi \,d\phi
\int\limits_0^{\infty}\frac{t^{m-n}dt}{(t/\rho_0+1)^{m+1}}=\\ 
=\rho_0^{m-n+1}B(m-n-p+1/2, p+1/2)B(m-n+1,n),
\end{multline*}
where $\rho_0=\rho(0,0)$, and $B$-functions stand for the values
of the first and the second integral in the product respectively
(cf. [GR], formulas 3.621(5) and 3.194(3); note that formula 4.368(3)
which could be used directly for computation of the integral in (3.9)
is inaccurate).
Expressing the $B$-functions in terms of $\Gamma$-functions we find that
(3.9) may be represented as $(-\lambda)^{-n}I_{mnp}$, where:  

$$
I_{mnp}=
(-1)^{m-n}
\frac{1}{4\pi^2}
\rho_0^{m-n+1}
\frac{\Gamma(m-n-p+1/2)\Gamma(p+1/2)\Gamma(n)}{\Gamma(m+1)}.
$$

Since $|\alpha|$ should be  positive, we have $m>n$.
On the other hand, as it was proved in ([AK])  
the order of every operator $[B,A;J|$ is not greater
than $|J|+2r$ and hence $|\alpha| \le|J|+2r$ which implies
$|J| \le 2n$. Another statement from ([AK]) gives the same estimate
on the length of the vector: $r \le 2n$. Therefore if we are interested
in the coefficient $a_n$ it is sufficient to consider $X_m$ with $m \le 4n$. 

Together with (3.3), (3.5) and (3.6) this implies:
\begin{multline}
a_n(x)=
\frac{1}{(n-1)!}\rho(u,v)^{-1}\cdot\\
\cdot
\sum_{m=n+1}^{4n}\sum_{k=0}^m (-1)^k \binom mk \Delta^k 
\Delta_0^{m-k} (\sum_{p=0}^{m-n}I_{mnp}
\frac{u^{2m-2n-2p}v^{2p}}{(2m-2n-2p)!(2p)!})|_{u=0,v=0}
\end{multline}

Computing explicitly $\Delta_0^{m-k}$ using (3.8) we may apply
the obtained differential operator to the last sum in (3.10).
Afterwards we change for conveniency one of the summation indices 
and  omit the vanishing terms.
Finally we obtain:
     
\begin{equation}
a_n(x)=\sum_{m=n+1}^{4n}\sum_{k=n+1}^m \sum_{s=0}^{k-n}C_{nksm}
\rho(u,v)^{k-n}
\Delta^k\left(u^{2k-2n-2s}v^{2s}\right)|_{u=0,v=0},
\end{equation}

where the constants $C_{nksm}$ are given by 
\begin{equation}
C_{nksm}=\frac{(-1)^n}{4\pi^2}\sum_{l=0}^{m-k}
\frac{ \Gamma(k+l-n-s+1/2)\Gamma(s+m-k-l+1/2)}
{k!\, l! \, (m-k-l)!\, (2k-2n-2s)!\, (2s)!}.
\end{equation}
This completes the proof of the Theorem 1.4. \qed

\section{Curvature coordinates and an invariant expression for $a_n(x)$
in a ``generic'' case}
The commutator method essentially involves non--invariant 
coordinate representation.
Nevertheless, using a special substitution it is possible to get the result 
in an invariant form, though in such a way that  it is valid for 
in some sense ``generic'' germs of curvature functions on $2$-manifolds.
\begin{definition}
Let $x$ be a point on the manifold $M$, and
$K(x)=K_0$ be the Gaussian curvature at this point.
Let $(z,w)$ be coordinates near the point $x$ and suppose that 
in a neighborhood of $x$ the map 
$$(z,w) \to (z_K,w_K)=(K-K_0, \Delta K - \Delta K_0),$$ 
where $\Delta K_0=(\Delta K)(x)$, has a nonvanishing Jacobian.
The  coordinates $(z_K,w_K)$ are called the 
{\it curvature coordinates} in the neighborhood of the point $x$.
\end{definition}
\begin{lemma}
The operator $\Delta_0=\Delta_x'$ may be written 
in curvature coordinates as
$$
\Delta_0=-(E\partial^2/\partial z_K^2 + 2F\partial^2/\partial z_K \partial w_K + G\partial^2/\partial w_K^2),
$$
where
$$
2E=2g^{11}=-\Delta(z_K^2)=2K\Delta K-\Delta(K^2);
$$
$$
2F=2g^{12}=-\Delta(z_K w_K)=K\Delta^2 K + (\Delta K)^2 - \Delta(K\Delta K); 
$$
$$
2G=2g^{22}=-\Delta(w_K^2)=2\Delta K \Delta^2K - \Delta (\Delta K)^2,
$$
and all the values are taken at the point $x=(0,0)$.
\end{lemma}
\begin{lemma} 
In the new coordinates
$$u=\sqrt{EG-F^2}z_K \quad v=Ew_K-Fz_K,$$
the Riemannian metric $\{g_{ij}\}$
takes the following form at the point $x=(0,0)$:  
\begin{equation}
ds^2=\frac{1}{E(EG-F^2)}(du^2+dv^2),
\end{equation}
where $E$, $F$, $G$ are as in Lemma {\em 4.2.} 
\end{lemma}

Lemmas 4.2. and 4.3. are proved by easy computations.  

The proof of the Theorem 1.4.  uses the conformal
form of the metric in fact only at the point $x=(0,0)$, 
therefore we may apply this theorem to (4.4),  
and expressing the result in curvature coordinates we obtain:

\begin{theorem}
Suppose that curvature coordinates  exist in some neighborhood of the 
point $x=(0,0)\in M$.
Then the  heat invariants $a_n(x)$ 
are given by the following formula  
\begin{multline}
a_n(x)=\sum_{m=n+1}^{4n}\sum_{k=n+1}^m \sum_{s=0}^{k-n}\sum_{p=0}^{2s}
(-1)^p \binom{2s}{p} C_{nksm}
\frac{E^{n-k+p} F^{2s-p}}
{(GE-F^2)^s} \cdot\\
\cdot \Delta^k\left((K-K_0)^{2k-2n-p}(\Delta K -
\Delta K_0)^p\right)|_{x=(0,0)},
\end{multline}
where $C_{nksm}$ are the same as in
formula {\em (3.12)}
\end{theorem}

\medskip

Note that all terms in this expression are invariant --- 
indeed, it is constructed from the powers of the Laplacian and the 
Gaussian curvature taken at the given point.

It is not always possible to introduce 
curvature coordinates --- for example, on a standard round 
sphere the Gaussian curvature is constant and therefore such coordinates 
are degenerate.
Nevertheless, they do exist "generically":
\begin{theorem}
In the space $\operatorname{J}^r$ of $r$-jets ($r\ge 3$) of germs of analytic
functions $K(u,v)$ at the point $x=(0,0)$,
the set of $r$-jets for which the ``curvature coordinates'' are degenerate
forms an algebraic set of codimension at least $1$.
\end{theorem}

\noindent{\bf Proof.}
We demand the curvature function $K(x)$, $x\in M$ be analytic near the 
point $x\in M$ in some local coordinates (and hence in any).
One may show that every germ of an analytic function maybe taken as
a germ of some curvature function. Indeed, if we write the metric locally
in conformal form (3.7) the Gaussian curvature is
given by 
\begin{equation}
K(u,v)=\frac{1}{2}\Delta(\log (\rho(u,v)).
\end{equation}
Clearly, for every analytic function $K(u,v)$ there exists a function
$\rho(u,v)$ satisfying the above PDE with initial conditions
$K(0,0)=K_0$ (by Cauchy--Kovalevskaya theorem --- see, for example,
[ES]).

Consider the space ${\mbox J}^r$ of $r$-jets ($r \ge 3)$ 
of germs of analytic functions
$K(u,v)$ at the point $x=(0,0)$.
Non--degeneracy of the ``curvature coordinates'' is equivalent to
the nonvanishing condition for the determinant of the Jacobi matrix 
$\operatorname{Jac}(K,\Delta K)$ at the
point $x$. The equation  
$\operatorname{Jac}(K,\Delta K)=0$ is an algebraic equation in
the space ${\mbox J}^r$, 
therefore (since the equality is obviously not identical) 
the set of its solutions is an algebraic set of codimension 
at least $1$ (see [Br]) in the space ${\mbox J}^r$. \qed

\medskip

\section* {Appendix}
As a further check of the formula (3.11) we programmed it 
(together with its slight modification for diagonalized metrics) 
using  the {\it Mathematica} software ([Wo]).

Here is the result of the work of the program for $a_1(x)$:
$$a_1(x)=\frac{1}{24\pi}\frac{\rho_u^2+\rho_v^2-\rho_{uu}\rho-\rho_{vv}\rho}
{\rho^3}$$
On the other hand, using the formula (4.8) we get that
$$K=\frac{\rho_u^2+\rho_v^2-\rho_{uu}\rho-\rho_{vv}\rho}
{2\rho^3},$$
and hence $$a_1(x)=\frac{K}{12\pi}=\frac{\tau}{24\pi},$$
where $\tau=2K$ is the scalar curvature at the point $x$.
This coincides with the well--known formula for $a_1(x)$ (cf. [G3]).

Let us emphasize that our aim was just to {\it test} the formula (3.11);
for {\it computational} purposes {\it Mathematica} is not an efficient tool 
for our method --- it does not manage to calculate even $a_2(x)$ for
a general metric in reasonable time (though in particular
cases of a ``simple'' conformal metric with 
$\rho(u,v)=1/(a_0+a_1u+a_2v)$, or a spherical
metric $ds^2=R^2(d\theta^2+\sin^2\theta d\phi^2)$ --- here we used a 
modification of (3.11) for diagonalized metrics, --- 
we found $a_2(x)$, and the result also agreed with the 
already known (see [G1],[G3])).
The main difficulty  is to compute high  powers of the Laplacian, which 
is a differential  operator with non--constant coefficients, and  
therefore the number of  terms in the expression for $\Delta^k$  
grows exponentially in $k$. However, we think  that programming the 
commutator method  in a language like C or in some  lower level 
language one may obtain  completely  explicit expressions for higher 
heat invariants and, hopefully, even for some unknown ones.

\bigskip

\bigskip
 
\bigskip

\centerline {\Large References}

\bigskip

[ABC] P.Amsterdamski, A.Berkin, and D.O'Connor, $b_8$ ``Hamidew''
coefficient for a scalar field, Class. Quant. Grav. vol.6, (1989), 1981-1991.

[AK] S.Agmon and Y.Kannai, On the asymptotic behavior of spectral 
functions and resolvent kernels of elliptic operators, Israel J. Math. 5
(1967), 1-30.

[Av] I.V.Avramidi, The covariant technique for the calculation of the 
heat kernel asymptotic expansion, Phys. Let. B 238 (1990), 92-97.

[Be] M.Berger, Geometry of the spectrum, Proc. Symp. Pure Math. 27 (1975),
129-152.

[Br] Th. Br\"ocker, Differential germs and catastrophes, London Math. Soc.
Lect. Note Series, Cambridge Univ. Press, (1975). 

[DNF] B.A.Dubrovin, S.P.Novikov and A.T.Fomenko, Modern geometry: methods
and applications, Nauka, (1986) (in Russian).

[ES] Yu.V. Egorov, M.A. Shubin, Partial Differential Equations I,
Encycl. Math. Sci. vol. 30, Springer--Verlag, (1992).

[F] S.A. Fulling ed., Heat Kernel Techniques and Quantum Gravity,
Discourses in Math. and its Appl., No. 4, Texas A\&M Univ., (1995).

[G1] P.Gilkey, The spectral geometry of a Riemannian manifold, J. Diff.
Geom. 10 (1975), 601-618.   

[G2]  P.Gilkey, The index theorem and the heat equation, Math. Lect. 
       Series,  Publish or Perish, (1974).

[G3] P.Gilkey, Heat equation asymptotics, Proc. Symp. Pure Math. 54 (1993),
317-326

[G4] P. Gilkey, Invariance theory, the  heat equation and the
Atiyah--Singer index theorem, Math. Lect. Series, Publish or Perish, (1984)

[GR] I.S. Gradshtein, I.M. Ryzhik, Table of integrals, series and
products, Academic Press,  (1980).

[Kan] Y.Kannai, On the asymptotic behavior of resolvent kernels, spectral
functions and eigenvalues of semi--elliptic systems, 
Ann. Sc. Norm. Sup. Pisa, Classe di Scienze  XXIII, (1969), 563-634.

[Sa] T.Sakai, On the eigenvalues of the Laplacian and curvature of Riemannian
manifold, T\^ohuku Math. J. 23 (1971), 585-603.

[vdV] A.E.M. van de Ven, Index--free heat kernel coefficients,
hep-th/9708152, (1997), 1--38.

[Wo] S.Wolfram, Mathematica: a system for doing mathematics by computer,
Addison--Wesley, (1991).

[Xu] C.Xu, Heat kernels and geometric invariants I, Bull. Sc. Math.
     117 (1993), 287-312.

\end{document}